\documentstyle[xypic,epic,amsthm,amsfonts,amsbsy,amssymb,amsgen,
amsmath,amsopn,verbatim,11pt]{article}
\input epsf
\input amssymb.sty

\begin{document} 

\theoremstyle{plain}
\newtheorem{thm}{Theorem}[section] 
\newtheorem{lem}[thm]{Lemma}
\newtheorem{prop}[thm]{Proposition}
\newtheorem{cor}[thm]{Corollary}
\newtheorem{defn}[thm]{Definition}
\newtheorem{problem}[thm]{Problem}

\theoremstyle{definition}
\newtheorem{rem}[thm]{Remark}
\newtheorem{egs}[thm]{Examples}
\newtheorem{eg}[thm]{Example}

\newcommand{\suredot}{\ }
%\newcounter{subs}[section]
%\renewcommand{\subsection}
%{\addtocounter{subs}{1}\noindent{\bf \thesection.\thesubs.\ }}

\addtolength{\textwidth}{3cm}%\addtolength{\hoffset}{1.5cm}
\addtolength{\textheight}{3.2cm} \addtolength{\voffset}{-1.6cm}

\author{}
\title{Grassmannian Trilogarithms}

\maketitle
\begin{center}
{\large Alexander B. Goncharov}
\end{center}
\begin{center}
{Department of Mathematics, Brown University, RI 02912}
\end{center}

\begin{center}
{\large Jianqiang Zhao}
\end{center}
\begin{center}
{Department of Mathematics, University of Pennsylvania, PA 19104}
\end{center}

{\bf Abstract}. In \cite{Gadv,Ggalois} and \cite{Gchow} two completely 
different constructions of single valued Grassmannian trilogarithms 
were given. The construction of \cite{Gchow} is very simple 
and provides Grassmannian $n$-logarithms for all $n$. 
However its motivic nature is hidden. The construction in 
\cite{Gadv,Ggalois} is very explicit and motivic, but its generalization 
for $n>4$ is not known. 
 
In this paper we will compute explicitly the Grassmannian trilogarithm
constructed in \cite{Gchow} and prove that it differs from the motivic 
Grassmannian trilogarithm by an explicitly given 
product of logarithms. We also derive some general results about the 
Grassmannian polylogarithms.

\section{Introduction}

\subsection{\bf The Grassmannian $n$-logarithm\suredot}\label{iso}
Let $V_m$ be an $m$-dimensional vector space over an arbitrary 
field $F$ with a given basis $e_0, \dots, e_{m-1}$.
Let $\{z_i\}$ be the coordinate system in $V_m$ dual to the 
basis $\{e_i\}$ and $\widehat{G}_p^q$ the Grassmannian of 
$p$-dimensional subspaces in generic position with respect to the
coordinate hyperplanes in $V_{p+q}$. The intersection of a 
hyperplane with the coordinate plane $z_i=0$  provides a map
$a_i: \widehat{G}_p^q \longrightarrow \widehat{G}_{p-1}^q$. 
The collection of the maps $\{a_i\}$ provides 
 a truncated semisimplicial variety over ${\Bbb Z}$
\begin{equation} \label{GG}
\text{$\diagramcompileto{diag1}
\text{$\widehat{G}_{n}^n$}\rto<1.25ex>_-{:}\rto<-1.2ex>^-{.}&
\ \text{$\widehat{G}_{n-1}^n$}\ \rto<1.25ex>_-{:}\rto<-1.2ex>^-{.}&
\ \cdots\ \rto<1.25ex>_-{:}\rto<-1.2ex>^-{.}&
\ \text{$\widehat{G}_{1}^n$}\
\enddiagram$}
\end{equation}
where $\widehat{G}_{n-k}^n$ sits in degree $2n-k$. 
Notice that $\widehat{G}_1^n = ({\Bbb G}_m)^{n}$. Indeed, it consists of the one dimensional subspaces in $V_{n+1}$ which do not lie in the hyperplanes $z_0=0$, ..., $z_n=0$. So 
$z_1/z_0, \dots, z_{n}/z_0$ are natural coordinates on $\widehat{G}_1^n$.

In \cite{Gchow} we constructed a collection of 
${\Bbb R}(n-1)=(2\pi i)^{n-1}{\Bbb R}$-valued differential 
$k$-forms ${\cal L}^{\rm G}_{k;n}$ on the complex Grassmannians 
$\widehat{G}_{n-k}^n({\Bbb C})$ satisfying the cocycle condition 
\begin{equation} \label{GG2}
d{\cal L}^{\rm G}_{k;n} = \sum_{i=0}^{2n-k-1} (-1)^{i} a_i^*{\cal L}^{\rm
G}_{k+1;n}\qquad -1\le k\le n-2
\end{equation}
and such that 
\begin{equation} \label{GG3}
d{\cal L}^{\rm G}_{n-1;n}(h_0,\dots,h_n) 
=-\pi_n\big( d \log (z_1/z_0) \wedge\dots\wedge d \log (z_{n}/z_0)\big)
\end{equation}
where $h_i=\{z_i=0\}$ and for any real numbers $a$ and $b$
$$\pi_n(a+bi)=
\cases
a \qquad &\text{$n$ odd} \\
bi&\text{$n$ even}. 
\endcases$$

A collection of forms as above is called (a {\it single valued})
Grassmannian $n$-logarithm. Its existence was conjectured in \cite{BMS}
(compare with \cite{GGL}). The function ${\cal L}^{\rm G}_{0;n}$ is called  
the Grassmannian $n$-logarithm function and usually denoted by
${\cal L}^{\rm G}_n$. 

{\bf Remark}. The condition (\ref{GG2}) for $k=-1$ is the 
$(2n+1)$-term functional equation for the Grassmannian $n$-logarithm function:
\begin{equation} \label{GG2**}
\sum_{i=0}^{2n} (-1)^{i} a_i^*{\cal L}^{\rm
G}_{0;n} = 0
\end{equation}

There are two other versions of the Grassmannian polylogarithms:  
the {\it real} Grassmannian polylogarithm function in \cite{GM}, 
on $\widehat G_{2n}^{2n}({\Bbb R})$, and the {\it multivalued complex analytic} 
Grassmannian polylogarithms on $\widehat{G}_{\bullet}^n({\Bbb C})$ 
in \cite{HM1, HM2}. In general the real Grassmannian polylogarithm is expected to live on 
$\widehat G_{\bullet}^{2n}({\Bbb R})$; it should be  responsible for the combinatorial 
Pontryagin classes, see 
[GGL] and [You]. 
We will not discuss them in our paper. The motivic construction of 
Grassmannian polylogarithms should in particular provide a coherent 
construction of all the three types of Grassmannian polylogs, 
as well as their $etale$, $p$-adic, etc. analogs. 

The coinvariants of the natural action of the group GL$(V_m)$ on the set of
all $n$-tuples of vectors in $V_n$ are called the {\it configurations} 
of $n$ vectors in $V_m$. The configuration spaces of $m$ 
vectors in two vector spaces of the same dimension are canonically 
isomorphic. So we only need to specify the dimension of the vector space 
when talking about configuration of vectors. We denote by
$C_m(V_n)$, or simply $C_m(n)$, the space of configurations of $m+1$ vectors 
in generic position in $V_n$. Then there are the following
{\it canonical} isomorphisms
\begin{equation} \label{duality}
\widehat{G}_p^q \cong C_{p+q-1}(V_p),\qquad
\widehat{G}_p^q \cong C_{p+q-1}(V_q).
\end{equation}
Namely, restricting the coordinate functions  $z_i$ to a subspace 
$W \in \widehat G^q_p$ we get a configuration of vectors 
$z_0,\dots, z_{p+q-1} \in W^*$. 
Projecting the vectors $e_i$ onto $V_{p+q}/W$ we get the second isomorphism.

\subsection{\bf Construction of the Grassmannian polylogarithms 
${\cal L}_{k;n}^{\rm G}$  \cite{Gchow}\suredot} 
First we need the following construction. Let $X$ be a variety
over ${\Bbb C}$ and $f_0,\dots,f_{n-1}$ be $n$ complex-valued functions on
$X({\Bbb C})$. We attach to the above data the following singular 
${\Bbb R}(n-1)$-valued differential $(n-1)$-form:
\begin{equation}\label{equ:rn}
r_{n}(f_1,\dots, f_n) :=-{\rm Alt}_{n}
\left\{\sum_{k\ge 0} c_{k,n}\log|f_1| \bigwedge_{j=2}^{2k+1}d\log|f_j| 
        \bigwedge_{j=2k+2}^n di\arg f_j \right\}
\end{equation}
where $c_{k,n}:={n\choose 2k+1}/n!$ and 
$${\rm Alt}_{n} F(x_1,\dots,x_n):=
\sum_{\sigma\in S_n}
{\rm sgn}(\sigma)F(x_{\sigma(1)},\dots,x_{\sigma(n)}).$$ 

The choice of the coefficients is dictated by the following property:
\begin{equation}\label{equ:drnForm}
d r_n(f_1,\dots,f_n)=-\pi_n
\left(d\log f_1\wedge\cdots\wedge d\log f_n\right)
\end{equation}

%When we choose linear coordinates $z_0,\dots, z_n$ in ${\Bbb CP}^n$
%such that the face $L_i$ of a simplex $L$ is given by $\{z_i=0\}$ and
%a hyperplane $H$ is given by $\sum z_i=0$ we see that
%$$r_n\big(\frac{z_1}{z_0},\dots,\frac{z_n}{z_0}\big)$$
%is a primitive of (see \cite[\S2]{Gchow})
%$$\pi_{n} \Omega_L=\pi_{n} \big(d\log\frac {z_1}{z_0},\dots,
 %       d\log\frac{z_n}{z_0}\big).$$

Let $l_0,\dots, l_{2n-k-1}$ be vectors in generic position in a 
complex vector space $V_{n-k}^*$. 
For $1 \leq i \le  2n-k-1$ set $f_i:= l_i/l_0$. They are 
$2n-k-1$ rational functions on ${\Bbb CP}^{n-k-1}$. 

\begin{defn} The Grassmannian $k$-form of weight $n$ on
$\widehat{G}_{n-k}^n({\Bbb C})$
is defined by
$$
{\cal L}_{k;n}^{\rm G}(l_0,\dots,l_{2n-k-1})=(2\pi i)^{k+1-n}
\int_{{\Bbb CP}^{n-k-1}} r_{2n-k-1}(f_1,\dots, f_{2n-k-1}).$$
\end{defn}
For the precise meaning of the right hand side see \cite{Gchow}
or \S\ref{setup} below. It was proved in \cite{Gchow} that this 
integral is convergent, so the definition makes sense. 

The Grassmannian $n$-logarithm function
${\cal L}_n^{\rm G}$ can be descended
onto the space of configurations of $2n$ points in 
$P(V^*_{n}) = {\Bbb CP}^{n-1}$ (see \cite{Gchow}) or, what is the same, 
the space of configurations of $2n$ hyperplanes in 
$P(V_n)$. Let $h_0$,\dots,$h_{2n-1}$ be $2n$ hyperplanes in 
${\Bbb CP}^{n-1}$. Choose rational functions $f_i$ such that
div$(f_i)=h_i-h_0$ for $1\le i\le 2n-1$. Then the Grassmannian 
$n$-logarithm function ${\cal L}_n^{\rm G}$ is defined by
$$
{\cal L}_n^{\rm G}(h_0,\dots,h_{2n-1})=(2\pi i)^{1-n}
\int_{{\Bbb CP}^{n-1}} r_{2n-1}(f_1,\dots, f_{2n-1}).$$

\subsection{\bf The Lie-motivic construction of the Grassmannian $n$-logarithms\suredot} 
A different construction of the the Grassmannian $n$-logarithms 
$L^{\rm G}_{\bullet;n}$ for $n=2,3$ 
was given  in \cite{Gadv, Ggalois} and for $n=4$ in \cite{Gaomoto},
see also \cite{GEK}. 

We will call these constructions {\it Lie-motivic} since 
they are obtained as a composition of a homomorphism from the Grassmannian 
complex (see \S\ref{review} below) to a motivic complex, understood 
as the weight $n$ part of the cochain complex of the motivic {\it Lie algebra},
followed by the canonical regulator map to the real Deligne complex.  

Let us explain in more details the notion of the Lie-motivic Grassmannian
polylogarithm function. It is expected that there is a natural variation 
of $n$-framed mixed Tate motives over the Grassmannian $\widehat G_n^n$ 
responsible for the motivic Grassmannian $n$-logarithm function in the 
following way. Taking the Hodge realization of this variation we get a 
variation of $n$-framed Hodge-Tate structures over the Grassmannian. 
Let ${\cal H}_n$ be the group of $n$-framed Hodge-Tate structures. 
Then ${\cal H}_{\bullet} = \oplus_{n \ge 0} {\cal H}_n$ has a natural 
Hopf algebra structure (see \cite{BGSV}). The coproduct on ${\cal H}_{\bullet}$ induces a graded Lie coalgebra structure on the quotient 
$$
{\cal L}({\cal H})_{\bullet}:= \frac{{\cal H}_{\bullet}}{{\cal H}_{>0}\cdot {\cal H}_{>0}}
$$
There are two natural period maps
\begin{equation} \label{BNM}
{P^H}: {\cal H}_{\bullet} \longrightarrow {\Bbb R}; 
\qquad  p^L : {\cal H}_{\bullet} \longrightarrow {\Bbb R}
\end{equation}
The first one is an algebra homomorphism, while the second kills the products:
$
p^L\Bigl({\cal H}_{>0}\cdot {\cal H}_{>0}\Bigr) = 0
$. 
Thus we get a canonical map $p^L: {\cal L}({\cal H})_{\bullet} \longrightarrow  {\Bbb R}$. 
Applying pointwise this map we get a function on the 
Grassmannian which we call the Lie-motivic Grassmannian polylogarithm $L^{\rm G}_n$. 

\subsection{\bf The comparison problem\suredot} Now a natural question arises:

\begin{problem}

a) What is the relation between the Grassmannian $n$-logarithms  
${\cal L}^{\rm G}_{\bullet;n}$ and $L^{\rm G}_{\bullet;n}$?
Do they coincide or not? 

b) Is it true that the Grassmannian $n$-logarithm 
${\cal L}^{\rm G}_{\bullet;n}$ admits a motivic construction? 

c) Is it true that the Grassmannian $n$-logarithm 
${\cal L}^{\rm G}_{\bullet;n}$ is Lie-motivic?
\end{problem}

By the very definitions one has ${\cal L}^{\rm G}_{n-1;n}=L^{\rm G}_{n-1;n}$. 

It was known from \cite{Gchow}, and it is already a nontrivial fact,  
that the Grassmannian dilogarithms of both types coincide. We will 
recover this result in s. 4.5 below. 

It was noticed by the first author during the preparation of \cite{Gchow}, 
and puzzled him very much, that the Grassmannian $n$-logarithms  
${\cal L}^{\rm G}_{\bullet;n}$ for $n \geq 4$ should be different from 
$L^{\rm G}_{\bullet;n}$. The reason is that 
${\cal L}^{\rm G}_{\bullet;n}$ satisfy some additional functional
equations which should not be true for $L^{\rm G}_{\bullet;n}$ for 
$n\geq 4$. Namely, projection along the subspace generated by 
$e_j$ provides a map 
$$
b_j: \widehat G_p^{q+1} \longrightarrow \widehat G_p^{q}.
$$
It was proved in \cite{Gchow} that ${\cal L}^{\rm G}_{p; q}$ satisfies
the property 
\begin{equation} \label{dual}
\sum_{j=0}^{2q -p}(-1)^j b_j^* {\cal L}^{\rm G}_{p; q} = 0.
\end{equation}
This property is valid for the motivic Grassmannian trilogarithm. 
However it should not be satisfied by $L^{\rm G}_{n-2;n}$ for $n\geq 4$, and because of this 
to construct the regulator map one  needs to extend the $L^{\rm G}$ 
to a {\it bi-Grassmannian} $n$-logarithm, see \cite{CHCL}. 

In this paper we compute explicitly the Grassmannian trilogarithm   
${\cal L}^{\rm G}_{k; 3}$ and show that it is different from 
$L^{\rm G}_{k; 3}$ for $k=0$ and $1$. The difference is explicitly computed
and has a motivic origin. Therefore the answer to 
part b) of the problem is positive for $n=3$. However it is not Lie-motivic, thus the answer 
to the question c) is negative.

\subsection{\bf Main result: computation of the Grassmannian trilogarithm 
${\cal L}^{\rm G}_3$\suredot} Recall that the classical polylogarithms 
are defined by
$$
Li_1(z) := -\log(1-z); \qquad Li_n(z):=\int_0^z Li_{n-1}(t)
\frac{dt}{t}, \quad n\ge 2.
$$
They admit the single-valued cousins (see \cite{Z}). For the dilogarithm
it is the Bloch-Wigner function
$$ \widehat {\cal L}_2(z) := i{\cal L}_2(z):= 
\pi_2\left( Li_2(z)\right)+ i \arg(1-z)\cdot \log|z|$$
and for the trilogarithm it is 
$$ {\cal L}_3(z)= {\mathfrak R}{\rm e} \big\{Li_{3}(z) - \log|z|
\cdot Li_2(z) \big\}- \frac13(\log|z|)^2\log|1-z|
$$ 
which was used in the proof of Zagier's conjecture on $\zeta(3)$ in 
\cite{Gadv} and \cite{Ggalois}. 

Choose a volume form $\omega\in\det V_n^*$ and set
$\Delta_{\omega}(v_1,\dots,v_n):=
\langle \omega, v_1\wedge\dots\wedge v_n\rangle \in F.$
We will usually omit the subscript $\omega$. 

The following result was proved in \cite{Gadv,Ggalois}, see also 
\cite{Gdenninger}.
\begin{equation}\label{mainL} 
{L}_3^{\rm G}(l_0,\dots, l_5) =  \frac{1}{90}{\rm Alt}_6 {\cal L}_3\left(
\frac{\Delta(l_0,l_1,l_3)\Delta(l_1,l_2,l_4)\Delta(l_2,l_0,l_5)}
{\Delta(l_0,l_1,l_4)\Delta(l_1,l_2,l_5)\Delta(l_2,l_0,l_3)}\right) 
\end{equation}

The next
theorem is proved in \S5:

\begin{thm}\label{main}
\begin{align*}
{\cal L}_3^{\rm G}(l_0,\dots, l_5) = &{L}_3^{\rm G}(l_0,\dots, l_5) \\
 - & \frac 19 {\rm Alt}_6 \Bigl(
\log |\Delta(l_0, l_1, l_2)| \log |\Delta(l_1, l_2, l_3)|
\log |\Delta(l_2, l_3, l_4)|\Bigr).
\end{align*}
\end{thm}

Our next goal is to show that {\it the function ${\cal L}_3^{\rm G}$ does not 
satisfy the most interesting functional equation valid for the
function $L_3^{\rm G}$}. Consider the following configuration of $6$ 
points on the projective plane, called the special configuration. 

\begin{center}
\hspace{4.0cm}
\epsffile{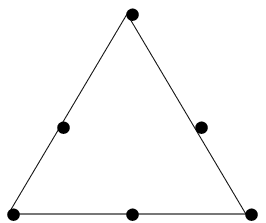}
\end{center}

The set of the special configurations can be naturally identified  with 
${\Bbb P}^1 \backslash \{0, \infty\}$ (see \cite{GPOL}). Namely $z \in F^*$ 
corresponds to the configuration $g_3(z)$ given by the columns of the 
following matrix 
$$\begin{pmatrix}
1 & 0& 0 & 1 & 0 & 1 \\ 
0 & 1& 0 & 1 & 1 & 0 \\ 
0 & 0& 1 & 0 & 1 & z
\end{pmatrix}$$
If a function is defined at the points of a set $X$, we can extend it 
by linearity to a homomorphism from the free abelian group generated by 
the points of $X$. We apply this construction to the functions $L_3^{\rm G}$
and ${\cal L}_3^{\rm G}$ and appropriate sets of the configurations of $6$ 
points in ${\Bbb P}^2$. Then according to \eqref{mainL} one has 
\begin{equation}\label{mainL1} 
L_3^{\rm G} \left( (l_0,\dots, l_5) - \frac{1}{90}{\rm Alt}_6 g_3\Bigl(  
\frac{\Delta(l_0,l_1,l_3)\Delta(l_1,l_2,l_4)\Delta(l_2,l_0,l_5)}
{\Delta(l_0,l_1,l_4)\Delta(l_1,l_2,l_5)\Delta(l_2,l_0,l_3)}\Bigr)\right)=0.
\end{equation}
However the function  ${\cal  L}_3^{\rm G}$ does not satisfy this functional 
equation. Indeed, we have the following result, which is proved in 
\cite{GPOL}, see also \cite{Gchow}. 

\begin{thm} The restriction of the Grassmannian trilogarithm function 
${\cal L}_3^{\rm G}$ to the special stratum coincides with the classical 
trilogarithm  function ${\cal L}_3$. More precisely,
${\cal L}_3^{\rm G}(g_3(z)) = {\cal L}_3(z)$. 
\end{thm}

Thus if the function ${\cal  L}_3^{\rm G}$ satisfies the functional equation  
\eqref{mainL1} it would coincide with $L_3^{\rm G}$ at the generic
configuration, which is not true according to Theorem \ref{main}. 

Notice that the function ${\cal L}_3^{\rm G}$ can be defined for an arbitrary 
configuration of points in ${\Bbb P}^2$ (\cite{Gchow}, \cite{GPOL}),
and one can show that it is continuous near the special stratum. 

{\bf Remark}. This shows that  in the \cite[\S 4]{Gaomoto} 
to define the motivic Lie coalgebra $\oplus G_{*}(F)$ 
of a field $F$ (cf. {\em loc. cit.}) one needs the functional equations 
for the Lie-motivic  Grassmannian $n$-logarithm $L^{\rm G}_n$, 
instead of its relative ${\cal L}^{\rm G}_n$.

{\bf Problem}. {\it Find explicit expression of the Lie-motivic Grassmannian
$n$-logarithm via the functions ${\cal L}^{\rm G}_n$}. 

\subsection{\bf Main result from the point of view of the Deligne cohomology} 
 
Recall that the Deligne cohomology $H^*_{{\cal D}}(X, {\Bbb R}(p))$ of a 
regular algebraic variety $X$ over ${\Bbb C}$ can be defined as the 
hypercohomology of $X$ with coefficients in the following complex of 
sheaves $\underline {\Bbb R}_{{\cal D}}(p)$ on $X(\Bbb C)$:
$$
\begin{array}{ccccccccccc} \label{del}
\Bigl({\cal D}_{X}^{0}&\hskip-.7ex\stackrel{d}{\longrightarrow}&\hskip-.7ex
{\cal D}_{X}^{1}&\hskip-.7ex\stackrel{d}{\longrightarrow}&\hskip-.7ex
\cdots&\hskip-.7ex\stackrel{d}{\longrightarrow}&\hskip-.7ex
{\cal D}^p_{X}&\hskip-.7ex\stackrel{d}{\longrightarrow}&\hskip-.7ex{\cal D}_{X}^{p+1}&\hskip-.7ex
\stackrel{d}{\longrightarrow}&\hskip-.7ex\cdots\Bigr) \otimes {\Bbb R}(p-1)\\
&\hskip-.7ex&&&&&{\ }\uparrow{\scriptstyle \pi_p} 
&&\hskip-.7ex{\ }\uparrow{\scriptstyle \pi_p}&&\\
&&&&&&\Omega^p_{X}
&\hskip-.7ex\stackrel{\partial}{\longrightarrow}
&\hskip-.7ex\Omega_{X}^{p+1}&\hskip-.7ex\stackrel{\partial}{\longrightarrow}&\hskip-.7ex
\cdots \phantom{\Bigr) \otimes {\Bbb R}(p-1)}
\end{array}
$$
Here ${\cal D}_{X}^{i}$ is the sheaf of $i$-currents on $X$. 
The group ${\cal D}_{X}^{0}$ is in degree $1$. To compute 
the 
Deligne cohomology of $X$ let us  replace $F^p\Omega^*$ by its Dolbeault resolution and denote by 
$\underline {\Bbb R}_{{\cal D}}(p)(X)$ the complex of the global sections of the complex of sheaves 
over $X(\Bbb C)$. Then 
$$
H^*_{{\cal D}}(X, {\Bbb R}(p)) = H^*(\underline {\Bbb R}_{{\cal D}}(p)(X))  
$$
Thus to calculate the 
Deligne cohomology $H^*_{{\cal D}}(\widehat G^n_{\bullet}, {\Bbb R}(n)) $ of the semisimplicial Grassmannian \eqref{GG} 
one needs to consider the cohomology of the total complex associated with the 
bicomplex of the shape
\begin{equation} \label{GG1}
\underline {\Bbb R}_{{\cal D}}(n)(\widehat G^n_n ) \longleftarrow \cdots
\longleftarrow \underline {\Bbb R}_{{\cal D}}(n)(\widehat G^n_1).
\end{equation}
A collection of differential forms on the simplicial Grassmannian 
\eqref{GG} satisfies the two conditions \eqref{GG2} and \eqref{GG3} 
if and only if it represents a $2n$-cocycle 
in the bicomplex \eqref{GG1}. 
By \cite[Lemma 2.3]{Gchow}, ${\cal L}^{\rm G}_{\bullet;n}$ satisfies
\eqref{GG2}. We will see that by definition, for $n=2$ and $3$, 
$L^{\rm G}_{\bullet;n}$ satisfies \eqref{GG2} and \eqref{GG3}, and moreover 
${\cal L}^{\rm G}_{n-1;n}=L^{\rm G}_{n-1;n}$.
Thus one can expect that the difference between the 
$6$-cocycles given by $L^{\rm G}_{\bullet;3}$ and ${\cal L}^{\rm G}_{\bullet;3}$ 
is a coboundary of a certain nice $5$-chain in the bicomplex \eqref{GG1}.
We shall prove that this is indeed the case and calculate explicitly the 
$5$-chain. It has a nonzero component only over $\widehat G_{2}^3$, 
and this component is given by the function
\begin{equation} \label{GG10}
(l_0,\dots, l_4) \longmapsto \frac 19 {\rm Alt}_6 \Bigl(
\log |\Delta(l_0, l_1, l_2)| \log |\Delta(l_1, l_2, l_3)|
\log |\Delta(l_2, l_3, l_4)| \Bigr).
\end{equation}
Denote by $C_3$ this $5$-chain 
in \eqref{GG1}. It is also of motivic nature: the function \eqref{GG10} 
is a composition of the map 
$$\aligned
\widehat{G}_2^3(F) &\longrightarrow S^3 F^* \\
 (l_0,\dots, l_4) &\longmapsto 
\frac 19 {\rm Alt}_6 \Bigl\{
\Delta(l_0, l_1, l_2)\cdot \Delta(l_1, l_2, l_3)
\cdot \Delta(l_2, l_3, l_4)\Bigr\},
\endaligned$$
which is defined for an arbitrary field $F$, with the logarithm homomorphism 
$$
S^3{\Bbb C}^* \longrightarrow {\Bbb R},\qquad 
x_1\cdot x_2\cdot x_3\longmapsto\log|x_1|\log|x_2|\log|x_3| 
$$ 
defined when $F = {\Bbb C}$. 

Denote by $\{{\cal L}^{\rm G}_{\bullet;3}\}$ and $\{{L}^{\rm G}_{\bullet;3}\}$
the $6$-cocycles in the bicomplex \eqref{GG1} provided by the 
collection of forms ${\cal L}^{\rm G}_{\bullet;3}$ and 
$L^{\rm G}_{\bullet;3}$. Let $D$ be the differential in the 
total complex associated with the bicomplex \eqref{GG1}. 

 \begin{thm} \label{TH2}
$$
\{{\cal L}^{\rm G}_{\bullet;n}\} - \{L^{\rm G}_{\bullet;n}\} = D(C_3).
$$ 
\end{thm} 

We expect a similar story for the Grassmannian $n$-logarithms in general: 
the forms ${\cal L}^{\rm G}_{\bullet;n}$ should have
a motivic nature in the following precise sense:  

1. One should have an explicitly given homomorphism ${\Bbb L}_n$ from the 
weight $n$ Grassmannian complex to the weight $n$ part 
$(\Lambda^*L(F)_\bullet^\vee, \Delta )_{(n)}$ 
of the cochain complex of the motivic Lie algebra $L(F)_{\bullet}$ 
of an arbitrary field $F$. (In fact ${\Bbb L}_n$ should be a part of 
the homomorphism from the bi-Grassmannian complex to the cochain complex 
of $L(F)_{\bullet}$). 

Composing this map with the regulator map  we get 
a cocycle in the bicomplex \eqref{GG1}, the  Lie-motivic 
Grassmannian $n$-logarithm ${L}^{\rm G}_{\bullet;n}$. 
For $n= 2,3,4$ this program has been implemented in \cite{Gadv},
\cite{Ggalois} and \cite{Gaomoto},  
but in general the homomorphism ${\Bbb L}_n$ is unknown. 
The story for for $n=2,3$ is recalled in section 4 below.  

2. We expect a natural $(2n-1)$-chain $C_n$ of motivic origin
in \eqref{GG1} such that 
$$
\{{\cal L}^{\rm G}_{\bullet;n}\} - \{L^{\rm G}_{\bullet;n}\} = D(C_n).
$$ 
Here $\{{\cal L}^{\rm G}_{\bullet;n}\}$ and $\{{L}^{\rm G}_{\bullet;n}\}$
are $2n$-cocycles in the bicomplex \eqref{GG1} provided by the 
forms ${\cal L}^{\rm G}_{\bullet;n}$ and 
$L^{\rm G}_{\bullet;n}$. 

Our desire to understand better the structure of Grassmannian 
polylogarithms was motivated by the following reasons:

i). The Grassmannian $n$-logarithm  can be used for an explicit 
construction of the class $c_n \in H^{2n}({\rm BGL}({\Bbb C})_{\bullet}, 
{\Bbb R}_{{\cal D}}(n))$ which provides Beilinson's regulator 
for ${\cal L}^{\rm G}$ (see \cite{GPOL}) and $L^{\rm G}$ (see
\cite{Gadv,Ggalois}). 

ii). Explicit calculation of the Grassmannian 
$n$-logarithm ${\cal L}^{\rm G}_{\bullet;n}$ 
should give a clue for construction  of the homomorphism 
${\Bbb L}_n$ as well as the chain $C_n$. 

{\bf Acknowledgment}.   
The first author gratefully acknowledges the support
of the NSF grant  DMS-9800998.

\section{Some properties of the differential forms $r_m$}

Here is another expression of the differential form $r_m$ which will be very
useful in applications. 
\begin{prop}\label{prop:rnholo} The differential $(m-1)$-form
$r_{m}(f_1,\dots, f_m)$ can be expressed as
\begin{equation}\label{i3}
{\rm Alt}_{m}\left\{\sum_{k=1}^m \frac {(-1)^{m-k-1}}{m!} 
  \log|f_1|\bigwedge_{j=2}^k d\log f_j
  \bigwedge_{j=k+1}^m d\,\overline{\log f_j}\right\}.
\end{equation}
\end{prop}
\begin{proof}
By definition
$$\aligned
r_{n+1} &(f_0,\cdots, f_n)=\frac {-1}{2^n(n+1)!} 
{\rm Alt}_{n+1}\left\{\sum_{j\ge 0}{n+1\choose 2j+1}
 \log|f_0| \cdot \right.\\
 \ &\ \left. \hskip2cm\bigwedge_{s=1}^{2j} \big(d\log f_s+d\,\overline{\log f_s}\big) 
 \bigwedge_{s=2j+1}^n \big(d\log f_s-d\,\overline{\log f_s}\big)\right\}\\
=&\frac {-1}{2^n(n+1)!} {\rm Alt}_{n+1}
\left\{\sum_{j\ge 0}{n+1\choose 2j+1}\sum_{k\ge 0} 
\sum_{l\ge 0} {n-2j\choose n-k-l}{2j\choose l} (-1)^{n-k-l} \right.\\
\ &\left.\ \hskip3cm\log|f_0|\bigwedge_{s=1}^k d\log f_s 
        \bigwedge_{s=k+1}^n d\,\overline{\log f_s}\right\}.
\endaligned$$
In the above, we've used the skew-symmetry property, for example, \\
${\rm Alt}_n\{d\,\overline{\log f_1}\wedge d\log f_2\cdots\}=
{\rm Alt}_n\{d\log f_1 \wedge d\,\overline{\log f_2}\cdots\}.$
The coefficients in the summation of index $k$ is obtained as follows:
for each $d\log f_s+d\,\overline{\log f_s}$ we can either
choose $d\log f_s$ or $\overline{\log f_s}$ but not both,
the same for $d\log f_s-d\,\overline{\log f_s}$.
For any appropriately fixed $l$, there are ${2j\choose l}$ 
ways to choose $\overline{\log f_s}$ from the former and
${n-2j\choose n-k-l}$ ways from the later. Once $\overline{\log f_s}$
are chosen, $d\log f_s$ are determined. We can now show that 
$$\sum_{0\le 2j\le n-1}{n+1\choose 2j+1} \sum_{l\ge 0} 
  {n-2j\choose n-k-l}{2j\choose l} (-1)^{n-k-l} =2^n (-1)^{n-k}$$
by comparing the coefficient of $x^{n-k}$ of the following polynomials
in $x$:
\begin{align*}
\ &\sum_{p=0}^n \sum_{0\le 2j\le n}{n+1\choose 2j+1} \sum_{l\ge 0} 
  {n-2j\choose p-l}{2j\choose l} (-1)^{p-l} x^p \\
=&\sum_{0\le 2j\le n}{n+1\choose 2j+1} (1-x)^{n-2j}(1+x)^{2j}\\
=&\frac 1{2(1+x)}\left[(1-x+1+x)^{n+1}-(1-x-1-x)^{n+1}\right]\\
=&\frac {2^n\left[1-(-x)^{n+1}\right]}{1+x} = 2^n \sum_{m=0}^n (-x)^m.
\end{align*}
The proposition follows at once.
\end{proof}
\begin{cor}\label{cor:rn} The $(2n-1)$-form $r_{2n}(f_1,\dots,f_{2n})$
can be expressed by
$${\rm Alt}_{2n} \pi_{2n}\left\{ \sum_{k=n+1}^{2n} \frac {2(-1)^{k-1}}{2n!}
        \log|f_1|\bigwedge_{j=2}^k d\log(f_j)
        \bigwedge_{j=k+1}^{2n} d\,\overline{\log(f_j)}\right\}
$$
and the $(2n-2)$-form $r_{2n-1}(f_1,\dots, f_{2n-1})$ is
\begin{multline*}
{\rm Alt}_{2n-1}{\em Re}\left\{\frac {(-1)^{n}}{(2n-1)!}
        \log|f_1|\bigwedge_{j=2}^{n} d\log(f_j)
        \bigwedge_{j=n+1}^{2n-1} d\,\overline{\log(f_j)}\right.\\
 +\sum_{k=n+1}^{2n-1} \frac {2(-1)^k}{(2n-1)!}
        \log|f_1|\bigwedge_{j=2}^k d\log(f_j)
        \left.\bigwedge_{j=k+1}^{2n-1} d\,\overline{\log(f_j)} \right\}.
\end{multline*}
\end{cor}
\begin{proof}
We can use symmetry to bring \eqref{i3} into the form in which
at least $\lceil (m-1)/2\rceil$ holomorphic $d\log$ appear together 
with at most $\lceil (m-2)/2\rceil$
anti-holomorphic $d\,\overline{\log}$.
\end{proof}

\begin{egs}\label{r3}
We will need the following special cases later:

\noindent(1) From Corollary \ref{cor:rn} $r_3(f_1,f_2,f_3)$ is equal to
$$\frac 16{\rm Re}\left({\rm Alt}_3
        \Bigl\{\log|f_1|d\log f_2 \wedge d\,\overline{\log f_3}-
        2\log|f_1|d\log f_2 \wedge d\log f_3 \Bigr\}\right).$$

\noindent(2) We have
$$\aligned
r_4(f_1,f_2,f_3,f_4)= &\frac1{12}\pi_4\left({\rm Alt}_4\Bigl\{
        \log|f_1| d\log f_2 \wedge d\log f_3\wedge d\,\overline{\log f_4}\right.\\
\ &\hskip2cm\left. -\log|f_1|d\log f_2\wedge d\log f_3\wedge 
        d\log f_4\Bigr\}\right).
\endaligned$$
\end{egs}

\section{Computation of the Grassmannian $1$-forms}\label{classG3log}

\subsection{\bf The setup\suredot}\label{setup}
For any $0\le k\le n-1$ we let $V_{n-k}$ be an 
$(n-k)$-dimensional complex vector space. 
Let $l_0,\dots, l_{2n-k-1}$ be vectors in generic position in 
the dual space $V_{n-k}^*$. Recall that
the Grassmannian $k$-form of weight $n$ is
\begin{equation}\label{equ:2nlogf}
{\cal L}_{k;n}^{\rm G}(l_0,\cdots,l_{2n-k-1})=(2\pi i)^{k+1-n}
\int_{{\Bbb CP}^{n-k-1}} r_{2n-k-1}(f_1,\dots, f_{2n-k-1}).
\end{equation} 
where $f_i(t)= l_i/l_0$, $1\le i\le 2n-k-1$, 
are $2n-k-1$ rational functions on $P(V_{n-k}) = {\Bbb CP}^{n-k-1}$. 
Our first goal is  to explain the meaning of  this integral. 
In the next subsection we give a recipe for its computation when $k=1$. 

Let $\pi: Z \longrightarrow Y$ be a smooth map of manifolds with compact fibers and 
$\omega$ be a distribution on $Z$. Then one can define $\pi_*\omega$ so that 
$\langle\pi_*\omega, \varphi\rangle=\langle\omega, \pi^*\varphi\rangle$
for any smooth test form $\varphi$ on $Y$. 

There is a canonical function on $V_{n-k}^* \times V_{n-k}$ whose value
at the point $(l,t)$ is $l(t)$. The expression 
$r_{2n-k-1}(f_1(t), \dots, f_{2n-k-1}(t))$ is a differential form with
logarithmic singularities on ${\Bbb CP}^{n-k-1}\times Y$ where
$$Y=\underbrace{V_{n-k}^* \times\cdots \times V_{n-k}^*}_{\text{$2n-k-1$
times}}.$$
It is proved in \cite{Gchow} that it has integrable singularities,
and thus provides a distribution on this manifold. 
The right hand side of \eqref{equ:2nlogf} is defined as
$(2\pi i)^{k+1-n}\cdot\pi_*\big(r_{2n-k-1}(f_1(t),\dots,f_{2n-k-1}(t))\big)$, 
where $\pi$ is the canonical projection along ${\Bbb CP}^{n-k-1}$. 
Write $d = d_t + d_a$ for the differential on $P(V_{n-k}) \times Y $,
where $d_t$ is the  $P(V_{n-k})$- and $d_a$ is the
$Y$-components of $d$.

Let $X$ and $Y$ be complex manifolds and $X$ is compact of complex 
dimension $d$. Let $\omega$ be a distribution on $X \times Y$. There 
is canonical projection $\pi: X \times Y \longrightarrow Y$. 
Denote by 
${\cal D}^{(p_1, q_1)}(X )$
the space of distributions of the Dolbeault type $(p_1, q_1)$. 
The space ${\cal D}(X \times Y)$ 
of distributions on $X \times Y$ admits a decomposition
$
{\cal D}(X \times Y) = \oplus {\cal D}^{(p_1,q_1; p_2, q_2)}(X \times Y), 
$
where $(p_1, q_1)$ (resp. $(p_2, q_2)$) is the type of the distribution
with respect to $X$ (resp. $Y$). If $\omega$ is of type
$(p_1, q_1; p_2, q_2)$ then $\pi_*\omega$ is of type 
$(p_1-d, q_1-d; p_2, q_2)$. In particular $\pi_*\omega=0$ 
if $p_1<d$ or $q_1<d$. 

Let us present $r_{2n-k-1}(f_1(t), \dots, f_{2n-k-1}(t))$ as a sum
of its Dolbeault components $\omega^{(p, 2n-k-2-p)}$. Then  
$$
\omega^{(p, 2n-k-2-p)} = 
\sum_{\alpha} \log|g^{\alpha}_0(t)| \bigwedge_{j=1}^p d\log g^{\alpha}_j(t) 
        \bigwedge_{j=p+1}^{2n-k-2} d\,\overline{\log g^{\alpha}_j(t)}
$$
where $g^{\alpha}_i(t)$ are some rational functions on ${\Bbb CP}^{n-k-1}$.
Therefore the integral of $\omega^{(p, 2n-k-2-p)}$ over
${\Bbb CP}^{n-k-1}$ is zero unless $n-k-1\le p\le n-1$. If $p=n-1$ 
the integral is calculated  as 
\begin{multline*}
\frac{{n-1\choose k}}{(2n-k-1)!}{\rm Alt}_{2n-k-1}\sum_\alpha\int_{{\Bbb CP}^{n-k-1}}
        \log|g^{\alpha}_0(t)|  \cdot \\
\bigwedge_{j=1}^{k}d_a \log g^{\alpha}_1(t) 
        \bigwedge_{j=k+1}^{n-1} d \log g^{\alpha}_j(t) 
        \bigwedge_{j=n}^{2n-k-2}d\,\overline {\log g^{\alpha}_j(t)}.
\end{multline*}
In this paper we are mostly interested in the case $k=1$.

\subsection{\bf The key formula\suredot} 
Our main task in this subsection is to calculate the Grassmannian $1$-form
${\cal L}^{\rm G}_{1;n+1}(l_0,\dots, l_{2n})$ of weight $n+1$. (We increase
the weight for ease of notation.)

Set $X:= V_n^* \times\cdots\times V_n^* $ and $Y:= V_n$. 
Let $(l_0,\dots, l_n; t)$ be a point of the variety
$X \times Y$. One has 
$$\omega^n(X \times Y) = \oplus_{a+b=n} \omega^a(X) \otimes 
        \omega^a(Y); \qquad \omega = \sum \omega^{(a; b)}.$$
We will now write $d = d_a+d_t$ where $d_t$ is the $V_n$-components of $d$. 
Let us compute the $(1;n-1)$ component of 
the following differential form:  
$$\sum_{i=0}^n(-1)^i d\log l_0(t) \wedge\dots\wedge d 
        \widehat{\log l_i(t)} \wedge \dots\wedge d\log l_n(t).$$
One can define the $SL(V_n)$-invariant Leray form
$$
\alpha_{n-1}(l_1(t),\dots,l_n(t)):= \sum_{i=1}^n(-1)^{i-1} 
        l_i(t) d_tl_1(t)\wedge\dots\wedge 
        \widehat{d_tl_i(t)}\wedge\dots\wedge d_tl_n(t).$$
Let $p: V_n \setminus \{l_i(t)=0\} \longrightarrow P(V_n)$ be the natural 
projection. Then one can check that the form 
$$
\frac{\alpha_{n-1}(l_1(t),\dots , l_{n}(t))}{l_1(t)\cdot\ldots\cdot l_{n}(t)}
$$
is lifted from $P(V_n)$, i.e. it is equal to $p^* \omega$ 
for some form $\omega$ on $P(V_n)$. 

\begin{prop} \label{prop:dlog=dgd}
\begin{align}
\ &\frac{1}{n!}{\rm Alt}_{(l_0,\dots,l_n)} \Bigl(d \log l_0(t) \wedge\dots 
\wedge d\log l_{n-1}(t)\Bigr)^{(1;n-1)} \notag\\
=&\frac 1{(n-1)!}{\rm Alt}_{(l_0,\dots,l_n)} \Bigl(d\log\Delta(l_0,\dots, l_{n-1})\wedge 
d_t \log l_1(t) \wedge\dots \wedge d_t\log l_{n-1}(t)\Bigr)\notag\\
=&\sum_{i=0}^{n}(-1)^{i} d \log \Delta(l_1,\dots, \widehat{l_i},\dots, l_n) \wedge 
\frac{\alpha_{n-1}(l_0(t),\dots, \widehat{l_i(t)},\dots, l_n(t))}
{l_0(t)\cdot\dots  \widehat{l_i(t)}\dots \cdot l_{n}(t)}.\label{LER1}
\end{align}
\end{prop}

\begin{rem} Because $\Delta(l_1,\dots, l_n)$ is a function on 
$V_n^* \times \cdots \times V_n^*$, one sees that 
$d_a\log \Delta(l_1,\dots, l_n) = d\log \Delta(l_1,\dots, l_n).$
\end{rem}

\begin{eg} 
In the simplest nontrivial case $n=2$ we get 
\begin{align*}
\ d\log l_0(t) &\wedge d\log l_1(t) - d\log l_0(t) \wedge d\log l_2(t)  
+ d\log l_1(t) \wedge d\log l_2(t) \\
=&d\log \Delta(l_0, l_1) \wedge \Bigl(d_t\log l_1(t) - d_t\log l_0(t)\Bigr) \\
-&d\log \Delta(l_0, l_2) \wedge \Bigl(d_t\log l_2(t) - d_t\log l_0(t)\Bigr) \\
+&d\log \Delta(l_1, l_2) \wedge \Bigl(d_t\log l_2(t) - d_t\log l_1(t)\Bigr).
\end{align*}
In coordinates it looks as follows:
\begin{align*}
 \ &d \log \frac{a_1 t_1+ a_2t_2}{c_1t_1 + c_2t_2}\wedge d \log 
        \frac{b_1 t_1+ b_2t_2}{c_1t_1 + c_2t_2}
=d\log(a_1 b_2-a_2 b_1)\wedge d_t\log \frac{b_1t_1+b_2t_2}{a_1t_1+a_2t_2}\\
-& d \log(a_1 c_2 - a_2 c_1)\wedge d_t\log 
        \frac{c_1t_1 + c_2t_2}{a_1t_1+a_2t_2} 
+ d\log(b_1c_2 - b_2c_1)\wedge d_t\log \frac{c_1t_1 +c_2t_2}{b_1t_1+b_2t_2}.
\end{align*}
\end{eg}
\begin{proof}
Choose a volume form $\omega \in \det V_n^*$. 
Denote by $\omega^{-1}\in \det V_n$ the dual volume form in $V^*_n$. 
Then for any vectors $l_1,\dots, l_n \in V_n^*$ we have 
$\Delta_{\omega^{-1}}(l_1,\dots,l_n ) \in F$. 
It is easy to check that
\begin{equation} \label{MYTO}
\alpha_{n-1}(l_1(t),\dots, l_n(t)) = 
\Delta_{\omega^{-1}}(l_1,\dots, l_n )\cdot i_E\omega
\end{equation}
where $E:= \sum t_i \partial_{t_i}$ is the Euler vector field in $V_n$. 
It follows from this that
\begin{multline*}
\Bigl(\sum_{i=0}^n(-1)^i\bigwedge_{0\le j\le n,j\ne i} 
        d \log l_j(t)\Bigr)^{(1;n-1)} \\
=\frac{\sum_{i=0}^n(-1)^i \Delta_{\omega^{-1}}(l_0,\dots, \widehat{l_i},\dots,l_n ) 
\cdot d_al_i(t)}{l_0(t)\cdots l_n(t) }\wedge i_E\omega.
\end{multline*}
Now let us calculate \eqref{LER1}. Using \eqref{MYTO} we get 
$$
d\log\Delta_{\omega^{-1}}(l_0,\dots,l_{n-1})
\frac{\alpha_{n-1}(l_0(t),\dots,l_{n-1}(t))}{l_0(t)\dots l_{n-1}(t)} 
=\frac{d\Delta_{\omega^{-1}}(l_0,\dots,l_{n-1})\wedge i_E\omega}{l_0(t)\dots l_{n-1}(t)}.
$$
So it remains to show that
$$
\sum_{i=0}^n(-1)^i l_i(t) d \Delta_{\omega^{-1}}
        (l_0,\dots, \widehat{l_i},\dots,l_n) 
=\sum_{i=0}^n(-1)^i \Delta_{\omega^{-1}}
        (l_0,\dots, \widehat{l_i},\dots,l_n)d_a l_i(t) 
$$
which follows by applying the differential $d_a$ to the identity
$$
\sum_{i=0}^n(-1)^i\Delta_{\omega^{-1}}
(l_0,\dots,\widehat{l_i},\dots,l_n)\cdot l_i(t)=0.
$$
We now can finish the proof by observing that
$$\alpha_{n-1}(l_0,\dots,l_{n-1})=\frac1{(n-1)!}{\rm Alt}_n
\Bigl\{l_0(t)d_tl_1(t)\wedge  \cdots \wedge d_tl_{n-1}(t)\Bigr\}.$$

\end{proof}

\begin{cor}\label{Grass1form}
 The Grassmannian $1$-form of weight $n+1$ is
\begin{multline*}
{\cal L}^{\rm G}_{1;n+1}(l_0,\dots, l_{2n})
=-\frac {(-2\pi i)^{1-n}}{(2n-1)!}\cdot
{\rm Alt}_{2n+1}\cdot \pi_{2n} \Bigl(d\log \Delta(l_1,\dots, l_n)\wedge \\
\int_{{\Bbb CP}^{n-1}} \log |l_0(t)| 
\bigwedge_{j=2}^n d_t\log l_j(t) \wedge \bigwedge_{j=n+1}^{2n-1} 
        d_t\,\overline{\log l_j(t)}\Bigr).
\end{multline*}
\end{cor}

\begin{proof} Let $f_i=l_i/l_0$ and 
$$b_n=\frac {2(-1)^{n}}{(2n)!}.$$
 From
$$r_{2n}(f_1, \dots, f_{2n})=
 \sum_{i=0}^{2n}(-1)^i r_{2n}(l_0, \dots, \widehat{l_i},\dots,l_{2n})$$
one has 
$$
{\cal L}^{\rm G}_{1;n+1}(l_0, \dots, l_{2n})=(2\pi i)^{1-n}
\int_{{\Bbb CP}^{n-1}}\frac{1}{(2n)!}{\rm Alt}_{2n+1} 
r_{2n}(l_0, \dots, l_{2n-1}).
$$
Using Corollary \ref{cor:rn} and observing that 
$\frac{1}{(2n)!}{\rm Alt}_{2n+1}{\rm Alt}_{2n}={\rm Alt}_{2n+1}$ 
we can rewrite the integrand as
\begin{align*}
\ &b_n\cdot{\rm Alt}_{2n+1}\pi_{2n} 
        \Bigl(\log|l_0(t)|\bigwedge_{j=1}^n d\log l_j(t)
                \bigwedge_{j=n+1}^{2n-1} d\,\overline{\log l_j(t)}\Bigr)\\
=&n\cdot b_n \cdot{\rm Alt}_{2n+1} \pi_{2n} \Bigl(\log |l_0(t)| 
d\log \Delta(l_1, \dots, l_n) \bigwedge_{j=2}^n d_t\log l_j(t)
         \bigwedge_{j=n+1}^{2n-1} d_t\,\overline{\log l_j(t)}\Bigr).
\end{align*}
\end{proof}

\section{The Grassmannian and polylogarithmic complexes: 
a review}\label{review}

\subsection{\bf The Grassmannian complex\suredot} 
Let $C_m(n)$ be the configurations of $m+1$ vectors in generic position
in $n$-dimensional vector space $V_n$ over $F$. Then there is a map 
$$
d': C_{m+1}(n+1) \longrightarrow C_m(n), \qquad 
(v_0,\dots,v_m) \longmapsto\sum_{i=0}^m 
(-1)^i(v_i|v_0,\dots,\widehat{v_i},\dots,v_m).
$$
Here $(v_i|v_0,\dots,\widehat{v_i},\dots,v_m)$ means the configuration
of $(v_0',\dots,\widehat{v_i'},\dots,v_m')$ in the space 
$V_n/\langle v_i\rangle$ where
$v_j'$ is the image of $v_j$ in $V_n/\langle v_i\rangle$.
It is straightforward to see that $(C_{*+n-1}(*),d')$ form a complex, called
the ($n$-th) Grassmannian complex. It is isomorphic to the 
complex $(C_{*+n-1}(n),d)$ where
$$
d':C_{m+1}(n) \longrightarrow C_m(n), \qquad 
(v_0,\dots,v_m) \longmapsto\sum_{i=0}^m (-1)^i(v_0,\dots,\widehat{v_i},\dots,v_m)
$$
by the duality $*: C_{m+n-1}(m)\rightarrow C_{m+n-1}(n)$
obtained by comparing the two isomorphisms in (\ref{duality}).

\subsection{\bf The polylogarithmic complexes\suredot} The polylogarithmic 
complex $(B(F;n)_\bullet, \delta)$ is a candidate to the weight 
$n$ motivic complex of the field $F$. 
It was defined in \cite{Gadv, Ggalois} for $n\leq 3$ as follows. 
The groups $B_n(F)$ are quotients $B_n(F):={\Bbb Z}[{\Bbb P}_F^1]/R_n(F)$, 
where the subgroups $R_n(F)$ reflect the known functional equations
for the $n$-logarithms for $n=1, 2,3$. 
For example, $B_1(F)=F^*\cong {\Bbb Z}[{\Bbb P}_F^1]/R_1(F)$ where
$R_1(F):=\big\langle\{x\} + \{y\} - \{xy\}:\ x,y\in F^*; 
\{0\};\{\infty\}\big\rangle$.  
Consider the homomorphisms
$${\setlength\arraycolsep{1pt}
\begin{array}{rccl}
\delta_n:& {\Bbb Z}[{\Bbb P}_F^1]&\longrightarrow& \cases
         F^* \wedge F^* \qquad&\text{if }n=2\\
         B_{2}(F)\otimes F^* &\text{if }n= 3
                \endcases \\
\ &\{x\} &\longmapsto & \cases
         (1-x)\wedge x \ \quad&\text{if }n=2\\
         \{x\}_{2}\otimes x &\text{if }n = 3
                \endcases\\
\ &\{0\},\{1\}, \{\infty\} &\longmapsto&  0
\end{array}}$$
where $\{x\}_{n}$ is the image of $\{x\}$ in
$B_{n}(F)$. Then $\delta_n(R_n(F)) =0$, so we get well defined homomorphisms 
$\delta: B_2(F) \longrightarrow \Lambda^2F^*$ and  $\delta: B_3(F) \longrightarrow B_2(F) \otimes F^*$.  We
get the polylogarithmic complexes
$$
B_2(F) \stackrel{\delta}{\longrightarrow} \Lambda^2 F^*, \qquad B_3(F)\stackrel{\delta}{\longrightarrow} 
B_2(F)\otimes F^*  \stackrel{\delta}{\longrightarrow}  \Lambda^3 F^*
$$
where $\delta(\{x\}_2 \otimes y) \longrightarrow (1-x) \wedge x \wedge y$. These complexes can
be thought of as the weight $2$ and $3$ parts of the standard 
cochain complex of the motivic Lie algebra $L(F)_{\bullet}$, see
\cite{Ggalois}. 

\subsection{\bf Homomorphisms from Grassmannian complexes to the
polylogarithmic ones\suredot}
There are  two commutative diagrams:
\begin{equation*}
\text{$\diagramcompileto{diagwt2}
C_3(2)\dto^-{\varphi_3(2)}\rto^-{d'}& C_2(1)\dto^-{\varphi_2(1)}\\
B_2(F)\rto^-\delta&\text{$\bigwedge$}^2 F^*
\enddiagram$}
\quad\text{ and }\quad\text{$\diagramcompileto{diagwt3}
C_5(3)\dto^-{\varphi_5(3)} \rto^-{d'}&C_4(2)\dto^-{\varphi_4(2)}\rto^-{d'}&
C_3(1)\dto^-{\varphi_3(1)}\\
B_3(F)\rto^-\delta&
B_2(F)\otimes F^*\rto^-\delta&\text{$\bigwedge$}^3 F^*
\enddiagram$}
\end{equation*}
where
$$
\varphi_2(1)(v_0,v_1,v_2)=
\frac 12{\rm Alt}_3\Bigl\{\Delta(v_0)\wedge\Delta(v_1)\Bigr\}$$
and $\varphi_3(2)(v_0,v_1,v_2,v_3)$ is given by
$\{r(v_0,v_1,v_2,v_3)\}_2$ which is the image of the cross-ratio 
\begin{equation} \label{CROSSR}
r(v_0,\dots,v_3)=\frac{\Delta(v_0,v_2)\Delta(v_1,v_3)}
        {\Delta(v_0,v_3)\Delta(v_1,v_2)}
\end{equation}
in $B_2(F)$.
For the second commutative diagram,
the map $\varphi_5(3)$ is the generalized cross-ratio
$$\varphi_5(3)(v_0,\dots,v_5)=\frac1{90}{\rm Alt}_6\Bigl\{\frac
{\Delta(v_0,v_1,v_3)\Delta(v_1,v_2,v_4)\Delta(v_2,v_0,v_5)}
{\Delta(v_0,v_1,v_4)\Delta(v_1,v_2,v_5)\Delta(v_2,v_0,v_3)}\Bigr\}_3\in
B_3(F)$$
and
\begin{align}
\varphi_4(2)(v_0,\dots,v_4)=&\frac 1{12} {\rm Alt}_5\Bigl\{r(v_0,v_2,v_3,v_4)_2
 \otimes \Delta(v_3,v_4)\Bigr\}\label{gf42}\\
\varphi_3(1)(v_0,\dots,v_3)=&-\frac 16{\rm Alt}_4\Bigl\{\Delta(v_0)
        \wedge\Delta(v_1)\wedge\Delta(v_2)\Bigr\}.\notag
\end{align}
\begin{rem} The correct proof of the second commutative diagram
was given in \cite{Gdenninger}. Notice that our $\varphi_5(3)$,
$\varphi_4(2)$ and $\varphi_3(1)$ are $1/6$ 
of the corresponding maps in \cite{Gdenninger}.
We made these changes in order to have
${\cal L}^{\rm G}_{n-1;n}=L^{\rm G}_{n-1;n}$.
\end{rem}

\subsection{\bf The regulator map on the polylogarithmic
complexes\suredot}\label{regmaps}
Let $X$ be a variety over ${\Bbb C}$ and $F:= {\Bbb C}(X)$. 
Let ${\cal A}_\eta^\bullet(X)$ is the de Rham complex of smooth 
forms at the generic point of $X$ over ${\Bbb C}$. 
Set $\alpha(f)=\log|f|d\log|1-f|-\log|1-f|d\log|f|$. Then 
there are the following commutative diagrams 
\begin{equation}\label{diag:wt=2}
\text{$\diagramcompileto{diagwt2}
B_2(F)\rto^-\delta\dto^{r_2(1)}&\text{$\bigwedge$}^2 F^*\dto^{r_2(2)}\\
{\cal A}_\eta^0(X) \rto & {\cal A}_\eta^1(X)
\enddiagram$}
\quad\text{ and }\quad\text{$\diagramcompileto{diagwt3}
B_3(F)\rto^-\delta\dto^{r_3(1)}&
B_2(F)\otimes F^*\rto^-\delta\dto^{r_3(2)}&\text{$\bigwedge$}^3 F^*\dto^{r_3(3)}\\
{\cal A}_\eta^0(X) \rto & {\cal A}_\eta^1(X)\rto &{\cal A}_\eta^2(X) 
\enddiagram$}
\end{equation}
where
$$\aligned
r_2(1) : \ & \{f\}_2 \longmapsto \widehat{{\cal L}}_2(f)=i{\cal L}_2(f)\\
r_2(2) : \ & g_0\wedge g_1\longmapsto r_2(g_0, g_1);\endaligned \
\aligned
r_3(1) : \ & \{f\}_3 \longmapsto {\cal L}_3(f)\\
r_3(2) : \ & \{f\}_2 \otimes g\longmapsto \widehat{{\cal L}}_2(f) d i\arg g 
       -\frac{\log|g|\alpha(f)}{3} \\
r_3(3) : \ & g_0\wedge g_1\wedge g_2\longmapsto r_3( g_0,g_1, g_2).
\endaligned
$$
Composing the maps $\varphi$ and $r$ we get 
the Lie-motivic Grassmannian polylogarithms $L^G_{\bullet;n}$ for $n=2,3$.

\subsection{\bf The Grassmannian dilogarithm\suredot} 
We deal with the left commutative diagram 
in \eqref{diag:wt=2}. It is easy to see that 
$$\aligned 
L^{\rm G}_{1;2}(l_0,l_1,l_2):= &r_2(2)\circ\varphi_2(1)(l_0, l_1, l_2)\\
=&\frac12{\rm Alt}_3\Bigl\{r_2(l_0,l_1)\Bigr\}
=r_2(f_1,f_2)
={\cal L}^{\rm G}_{1;2}(l_0, l_1, l_2).
\endaligned$$
Therefore the difference between $
{\cal L}^{\rm G}_2(l_0,\dots,l_3)
$ and 
$$
L^{\rm G}_{0;2}(l_0,\dots,l_3):= r_2(1)\circ\varphi_3(2)(l_0,\dots,l_3) = 
\widehat{\cal L}_2(r(l_0,\dots,l_3))
$$
is a constant. Since it is skewinvariant with respect to the permutations of the vectors $l_i$, it is zero, i.e. $L^{\rm G}_{0;2} = {\cal L}^{\rm G}_2$.

\section{Proof of theorems \ref{main} and \ref{TH2}}

By the very definition
$$\aligned 
L^{\rm G}_{2;3}(l_0,l_1,l_2,l_3):=&r_3(3)\circ\varphi_3(1)(l_0, l_1,l_2, l_3)\\
=&-\frac 16{\rm Alt}_4\Bigl\{r_3(l_0,l_1,l_2)\Bigr\}=r_3(f_1,f_2,f_3)
=&{\cal L}^{\rm G}_{2;3}(l_0,l_1,l_2,l_3).
\endaligned$$

We now want to compare
$$
L^{\rm G}_{1;3}(l_0,\dots,l_4):= r_3(2)\circ \varphi_4(2)(l_0, \dots, l_4)
$$
and 
\begin{multline*}
{\cal L}^{\rm G}_{1;3}(l_0,\dots,l_4)=(2\pi i)^{-1}
\int_{{\Bbb CP}^1} r_4(f_1, f_2, f_3, f_4) \\
=\frac{(2\pi i)^{-1}}{4!} 
{\rm Alt}_{(l_0,\dots, l_4)}\int_{{\Bbb CP}^1} r_4(l_0, l_1, l_2, l_3) 
\end{multline*}
Notice that the $(1,2)$-component of $r_4(l_0, l_1, l_2, l_3)$ 
is a $(1,2)$-form on $X \times V_2$, not on $X \times P(V_2)$.
However after the alternation we get a $(1,2)$-form on $X \times P(V_2)$. 

We will write $\Delta(a,b):=\Delta(l_a,l_b)$ for the rest of the paper.

\begin{prop}\label{grassLog1;3}
\begin{multline*}
{\cal L}^{\rm G}_{1;3}(l_0,\dots,l_4)={\rm Alt}_5\left\{\Bigl(
        \frac1{12}\widehat{{\cal L}}_2(r(l_0,l_1,l_2,l_4))di\arg \Delta (1,4)\right.\\
        \ \left.-\frac13  \log|\Delta(0,1)|\log|\Delta(1,4)|d \log|\Delta(2,4)|\Bigr)\right\}.
\end{multline*}
\end{prop}

\begin{proof} 
By Corollary \ref{Grass1form}
\begin{align}
6\cdot 2\pi i &{\cal L}_{1;3}^{\rm G}(l_0, l_1, l_2, l_3, l_4 )\notag\\\
=&{\rm Alt}_5 \pi_4\left\{ d\log \Delta(1,2)
 \int_{{\Bbb CP}^1}\log|l_0| d \log(l_2)\wedge d\,\overline{\log (l_3})\right\}\notag\\
=&-{\rm Alt}_5 \pi_4\left\{ d\log \Delta(1,4)
 \int_{{\Bbb CP}^1}\log|l_0| d \log(l_1)\wedge d\,\overline{\log (l_2})\right\}\notag\\
=&2{\rm Alt}_5 \left\{ d\log|\Delta(1,4)|
\int_{{\Bbb CP}^1}\log|l_0| d \log|l_1| \wedge di\arg l_2\right\}\label{MMM}\\
-&2{\rm Alt}_5  \left\{di\arg\Delta(1,4)\int_{{\Bbb CP}^1}
\log|l_0| d \log|l_1|\wedge d\log|l_2|\right\}\label{MMM4}
\end{align}
To get \eqref{MMM} and \eqref{MMM4} we use the following observations. 
Let $f$ and $g$ be holomorphic functions on a complex curve $X$ and 
$\varphi$ is a real valued function. Then since 
$\int_X \varphi d\log f\wedge d\log g =0$ we have, 
taking the real and imaginary parts respectively, 
$$\aligned
\int_X \varphi d\log|f| \wedge d\log |g| =
        &\int_X \varphi d\arg f \wedge d\arg g,\\ 
\int_X \varphi d\log|f| \wedge d\arg g =
        &\int_X \varphi d\log |g| \wedge d\arg f.
\endaligned
$$

One can easily show that (by using Examples \ref{r3}(1))
\begin{multline} \label{MMM10}
2\pi i {\cal L}_2^{\rm G}(l_0, l_1, l_2, l_3)=
\int_{{\Bbb CP}^1} r_3(f_1, f_2, f_3)=2\int_{{\Bbb CP}^1}
\log|f_1|d \log|f_2|\wedge d\log|f_3| \\
=-\frac{2}{3!}{\rm Alt}_4 
        \int_{{\Bbb CP}^1}\log|l_0|d \log|l_1|\wedge d\log|l_2|.
\end{multline}
This is a special case of equation \eqref{Grasslogn} in the Appendix. 

Writing \eqref{MMM4} as 
$$-{\rm Alt}_5\left\{di\arg\Delta(1,4)\int_{{\Bbb CP}^1}
        \Bigl({\rm Alt}_{l_1,l_4}\Bigl\{\log|l_0|
                d\log|l_1|\wedge d\log|l_2|\Bigr\}\Bigr)\right\}$$
and subtracting from this 
$$
{\rm Alt}_5\left\{di\arg\Delta(1,4)\int_{{\Bbb CP}^1}\Bigl( 
        {\rm Alt}_{l_1,l_2}\Bigl\{\log|l_0| d \log|l_1|
                \wedge d\log|l_4|\Bigr\} \Bigr)\right\}$$
which is zero (to check this use the skewsymmetry with respect 
the alternations $\{2,3\}$ and $\{0,3\}$), we see, using 
\eqref{MMM10}, that \eqref{MMM4} is equal to 
\begin{multline*}
{\rm Alt}_5\Bigl\{\pi i{\cal L}_2^{\rm G}(l_0,l_1,l_2,l_4) 
di\arg\Delta(1,4)\Bigr\}\\
=\pi i{\rm Alt}_5\Bigl\{\widehat{{\cal L}}_2(r(l_0,l_1,l_2,l_4))
di\arg\Delta(1,4)\Bigr\}
\end{multline*}
by \eqref{L2}.
To calculate \eqref{MMM} we compute, in two different ways, the expression
\begin{equation} \label{Nov9}
{\rm Alt}_5\left\{  d\log|\Delta(1,4)| \int_{{\Bbb CP}^1} d\log|l_1| \wedge   
d {\cal L}_2\Bigl(  \frac{\Delta(2,4)l_0}{\Delta(0,2)l_4}\Bigr)\right\}.
\end{equation} 

1. The integral over ${\Bbb CP}^1$, and hence the whole expression, is zero because
$$
d\log|l_1| \wedge   
d {\cal L}_2\Bigl(\frac{\Delta(2,4)l_0}{\Delta(0,2)l_4}\Bigr)
 = d\left\{ \log|l_1| \wedge   
d {\cal L}_2\Bigl(\frac{\Delta(2,4)l_0}{\Delta(0,2)l_4}\Bigr)\right\}
$$
where both parts are understood as currents. Notice that $\log |z|$ and 
${\cal L}_1(z)$ have integrable singularities and thus provide currents on ${\Bbb CP}^1$. 

2. Using formulas $\Delta(2,4)l_0(t)-\Delta(0,4)l_2(t)=\Delta(0,2)l_4(t)$ 
and 
$$ d{\cal L}_2(f)=\log|f|d\arg (1-f) -\log|1-f|d\arg f $$
we see that $d {\cal
L}_2\Bigl(\frac{\Delta(2,4)l_0}{\Delta(0,2)l_4}\Bigr)$ is equal to
\begin{multline*}
\log \left| \frac{\Delta(2,4)l_0}{\Delta(0,2)l_4}\right| d\arg
        \Bigl(\frac{\Delta(0,4)l_2}{\Delta(0,2)l_4}\Bigr)
-\log \left| \frac{\Delta(0,4)l_2}{\Delta(0,2)l_4}\right| d\arg 
        \Bigl(\frac{\Delta(2,4)l_0}{\Delta(0,2)l_4}\Bigr)\\
=\log \left| \frac{\Delta(2,4)l_0}{\Delta(0,2)l_4}\right| 
        d \arg \frac{l_2}{l_4}
-\log \left| \frac{\Delta(0,4)l_2}{\Delta(0,2)l_4}\right|
        d \arg \frac{l_0}{l_4}.
\end{multline*}
Since this expression is skewsymmetric with respect to the 
transposition $\{2,0\}$ exchanging the indices $2$ and $0$ 
we can write \eqref{Nov9} as 
\begin{align}
0=&2 {\rm Alt}_5\left\{ d\log|\Delta(1,4)| \int_{{\Bbb CP}^1} d\log|l_1| \wedge   
        \log\left|\frac{\Delta(2,4)l_0}{\Delta(0,2)l_4}\right| 
                d i\arg \frac{l_2}{l_4}\right\}\notag\\
=&2 {\rm Alt}_5\left\{d\log|\Delta(1,4)|\int_{{\Bbb CP}^1} 
        \log\left|\frac{l_0}{l_4}\right|d\log|l_1| 
                \wedge d i\arg \frac{l_2}{l_4}\right\} \label{03kill}\\
-&2 {\rm Alt}_5\left\{\log|\Delta(0,2)| d\log|\Delta(1,4)| \int_{{\Bbb CP}^1} 
        d\log|l_1|\wedge d i\arg \frac{l_2}{l_4}\right\}.  \label{MMM3}
\end{align}
We got last line by using transposition $\{0,3\}$.
Now \eqref{03kill} is exactly \eqref{MMM} since the other 
three possible terms vanish due to alternation $\{0,3\}$ (or $\{2,3\}$). 
Therefore \eqref{MMM} is equal to 
\begin{align}
\ &2{\rm Alt}_5\left\{\log|\Delta(0,2)| d\log|\Delta(1,4)| \int_{{\Bbb CP}^1} d\log|l_1| \wedge  
d i\arg \frac{l_2}{l_4} \right\} \notag \\
=&-2{\rm Alt}_5\left\{\log|\Delta(0,2)| d\log|\Delta(1,4)| \int_{{\Bbb CP}^1} \log|l_1|
\wedge d\big(d i\arg \frac{l_2}{l_4}\big) \right\} \notag \\
=&-4\pi i{\rm Alt}_5\left\{\log|\Delta(0,2)| d\log |\Delta(1,4)|
\log \left| \frac{\Delta(1,2)}{\Delta(1,4)}\right| \right\}. \label{lastexp}
\end{align}
We got this line by noting that
$d (d i\arg f) = 2\pi i \delta(f) $. Notice that
$$
4 \pi i{\rm Alt}_5\Bigl\{ \log|\Delta(0,2)| d\log|\Delta(1,4)|
\log|{\Delta(1,4)}| \Bigr\}=0
$$
since the expression is unchanged under the transposition $\{0,2\}$,
we see that \eqref{lastexp} equals to
$$
-4\pi i\cdot{\rm Alt}_5\Bigl( \log|\Delta(0,1)| d\log|\Delta(2,4)|\log|\Delta(1,4)|\Bigr)
$$
by transposition $\{1,2\}$ followed by $\{2,4\}$.
This finishes the proof of the proposition. 
\end{proof}

To calculate $L^{\rm G}_{1;3}(l_0,\dots, l_4)$ we need
\begin{prop}\label{logAlpha}
${\rm Alt}_5\Bigl(\log|\Delta(1,4)| \alpha(r(l_0, l_1, l_2,l_4))\Bigr)$
is equal to
\begin{multline*} 
4d{\rm Alt}_5\Bigl\{\log|\Delta(2,4)|\log|\Delta(1,4)|\log|\Delta(0,2)|\Bigr\} \\
+12{\rm Alt}_5\Bigl\{\log|\Delta(1,4)|\log|\Delta(0,1)|d\log|\Delta(2,4)|\Bigr\}.
\end{multline*}
\end{prop}

\begin{proof} Here is the algebraic reason behind this lemma. 
There is the following exact sequence of ${\Bbb Q}$-vector spaces
$$
F_{{\Bbb Q}}^* \otimes \Lambda^2 F_{{\Bbb Q}}^* \stackrel{\kappa_1}{\longrightarrow }
S^2 F_{{\Bbb Q}}^* \otimes F_{{\Bbb Q}}^* \stackrel{\kappa_2}{\longrightarrow } S^3F_{{\Bbb Q}}^*
$$
$$
\kappa_1: a \otimes b \wedge c \longmapsto a \cdot b \otimes c 
- a \cdot c \otimes b, \qquad \kappa_2: 
a\cdot b \otimes c \longmapsto a\cdot b\cdot c.
$$
It is a special case of the Koszul complex. 
The map $\kappa_2$ admits a natural splitting 
$$
\kappa'_2: a\cdot b\cdot c 
\longmapsto \frac13 \Bigl(a\cdot b \otimes  c  + a\cdot c \otimes  b  
        + b\cdot c \otimes  a  \Bigr).
$$
If $F={\Bbb C}(X)$ then there is a map
$$
S^2F^* \otimes F^* \longrightarrow {\cal A}^1(X); \qquad f_1 \cdot f_2 \otimes f_3 \longmapsto 
\log |f_1| \log |f_2| d \log |f_3|
$$

Now the proposition is an immediate corollary of the following lemma.
\end{proof}

\begin{lem}
\begin{align} 
-\kappa_1 {\rm Alt}_5&\Bigl\{\Delta(1,4) \otimes (1- r(l_0, l_1, l_2, l_4)) \wedge 
r(l_0, l_1, l_2, l_4)\Bigr\} \label{LEMMA3}\\
=&12\kappa'_2\left({\rm Alt}_5\Bigl\{\Delta(2,4) \cdot \Delta(1,4) \cdot 
        \Delta(0,2) \Bigr\}\right) \notag\\
+&12{\rm Alt}_5\Bigl\{\Delta(1,4) \cdot \Delta(0,1) \otimes \Delta(2,4) \Bigr\}.\notag
\end{align}
\end{lem}

\begin{proof} Let us show that that \eqref{LEMMA3} equals to 
\begin{multline} \label{LEM7}
4{\rm Alt}_5\Bigl\{\Delta(1,4)\cdot \Delta(0,2)  \otimes \Delta(0,1)-
        \Delta(1,4)\cdot \Delta(0,1)  \otimes \Delta(0,2)\\
+\Delta(1,4)\cdot \Delta(0,1)\otimes\Delta(1,2)\Bigr\}
\end{multline}
Indeed, using \eqref{CROSSR} we get 
$$
(1- r(l_0, l_1, l_2, l_4)) \wedge r(l_0, l_1, l_2, l_4) 
= \frac12{\rm Alt}_{(l_0, l_1, l_2, l_4)}\Bigl\{\Delta(0,1) \wedge \Delta(0,2)\Bigr\}
$$
Using this we write \eqref{LEMMA3} as a sum of the following $12$ terms:
$$\aligned
{\rm Alt}_5\Bigl(
-&\Delta(1,4)\cdot\Delta(0,1)\otimes\Delta(0,2)+\Delta(1,4)\cdot\Delta(0,2)\otimes\Delta(0,1)\\
+&\Delta(1,4)\cdot\Delta(0,1)\otimes\Delta(1,2)-\Delta(1,4)\cdot\Delta(1,2)\otimes\Delta(0,1)\\
-&\Delta(1,4)\cdot\Delta(0,2)\otimes\Delta(1,2)+\Delta(1,4)\cdot\Delta(1,2)\otimes\Delta(0,2)\\
-&\Delta(1,4)\cdot\Delta(0,2)\otimes\Delta(0,4)+\Delta(1,4)\cdot\Delta(0,4)\otimes\Delta(0,2)\\
+&\Delta(1,4)\cdot\Delta(0,2)\otimes\Delta(2,4)-\Delta(1,4)\cdot\Delta(2,4)\otimes\Delta(0,2)\\
-&\Delta(1,4)\cdot\Delta(0,4)\otimes\Delta(2,4)+\Delta(1,4)\cdot\Delta(2,4)\otimes\Delta(0,4)\Bigr)
\endaligned
$$
Notice that {\em a priori} \eqref{LEMMA3} is a sum of 24 terms
corresponding to the 24 terms in \eqref{CROSSR}. However 12 of 
them disappear after the alternation. For instance, 
${\rm Alt}_5(\Delta(1,4)\cdot\Delta(0,1)\otimes \Delta(0,4))=0$ since the 
involution $\{2,3\}$ does not change the expression. 

Computing the sign of the appropriate permutation, we see that this 
sum is equal to \eqref{LEM7}. Indeed, the terms $2$, $5$, $7$ and $9$
provide the first summand, the terms $1$, $6$, $8$, $10$ the second 
summand, and the rest the third summand. 

The third term in \eqref{LEM7} vanishes by using 
the involution $\{0,4\}$. This involution also
bring the first summand in \eqref{LEM7} into the following form
$4{\rm Alt}_5\Bigl\{
\Delta(1,4) \cdot \Delta(0,1) \otimes \Delta(2,4)\Bigr\}.$ 

The second term in \eqref{LEM7} contributes  
$$
8{\rm Alt}_5\Bigl\{\Delta(1,4) \cdot \Delta(0,1) \otimes \Delta(2,4)\Bigr\}+ 
12\kappa'_2\left(
{\rm Alt}_5\Bigl\{\Delta(2,4)\cdot\Delta(1,4)\cdot\Delta(0,2)\Bigr\}\right).
$$
The lemma, and hence Proposition \ref{logAlpha}, are proved. 
\end{proof}

Notice that $\{r(l_1,l_2,l_3,l_4)\}_2$ is skewsymmetric with 
respect to the permutations of $l_i$'s. So applying involution 
$\{1,3\}$ to \eqref{gf42} one gets
$$
\varphi_4(2)(l_0,\dots,l_4)=
\frac 1{12}{\rm Alt}_5\Bigl\{\{r(l_0,l_1,l_2,l_4)\}_2 \otimes \Delta
(l_1,l_4)\Bigr\}. 
$$
Therefore using Propositions \ref{grassLog1;3}
and the formula for $r_3(2)$ at the end of the \S\ref{regmaps}, we get 
\begin{multline}\label{keyequ}
{\cal L}^{\rm G}_{1;3}(l_0,\dots,l_4)-L^{\rm G}_{1;3}(l_0,\dots,l_4)
={\cal L}^{\rm G}_{1;3}(l_0,\dots,l_4)-
r_3(2)\circ \varphi_4(2)(l_0,\dots, l_4)\\
=\frac19 d\left({\rm Alt}_5\Bigl\{\log|\Delta(2,4)|
        \log|\Delta(1,4)|\log|\Delta(0,2)|\Bigr\}\right).
\end{multline}

Thus using this and formula (\ref{GG2}) in the case $k=0$, $n=3$ we conclude that 
$$
{\cal L}^{\rm G}_{0;3}(l_0,\dots,l_5) - L^{\rm G}_{0;3}(l_0,\dots,l_5) +
\frac19 {\rm Alt}_6\Bigl\{\log|\Delta(5,2,4)|
        \log|\Delta(5,1,4)|\log|\Delta(5,0,2)|\Bigr\}
$$
is a constant (notice the change of the sign before $1/9$). Since it is skewsymmetric with respect to the permutations of the vectors $l_0, ..., l_5$, it must be zero. 
Finally, notice that 
$$
{\rm Alt}_6\Bigl\{\log|\Delta(5,2,4)|
        \log|\Delta(5,1,4)|\log|\Delta(5,0,2)|\Bigr\} = 
$$
$$
{\rm Alt}_6\Bigl\{\log|\Delta(0,1,2)|
        \log|\Delta(1,2,3)|\log|\Delta(2,3,4)|\Bigr\}
$$
Theorem \ref{main} and \ref{TH2} are proved.

\section{Appendix: formulas for the Grassmannian $n$-logarithm function}
To simplify the Grassmannian $n$-logarithm we need

\begin{lem}\label{lem:xjyj}
Let $X$ be an $n$-dimensional complex manifold. Let
$f_1, \cdots, f_{2n}$ be any $2n$ rational functions 
on $X$. Then

\noindent{\sl{\em (I)}} For every $0 \le j \le n-1$ 
\begin{equation} \label{equ:xj} 
{\rm Alt}_{2n}\Bigl\{ \bigwedge_{k=1}^{2j+1} d\log|f_j|
        \bigwedge_{k=2j+2}^{2n} d\arg f_k \Bigr\} =0. 
\end{equation}
{\sl{\em (II)}} For every $0 \le j \le n-1$ 
\begin{equation}\label{equ:yj} 
{\rm Alt}_{2n}\Bigl\{\bigwedge_{k=1}^{2j} d\log|f_j|
        \bigwedge_{k=2j+1}^{2n} d\arg f_k \Bigr\}
=b_{j,n} d\log|f_1| \wedge \cdots \wedge d\log|f_{2n}|
\end{equation} 
where $b_{j,n}= (2n)!{n\choose j}/{2n\choose 2j}$.
\end{lem}

\begin{proof}
(I) Because $\dim X=n$, for any $0\le i\le n-1$ one has 
\begin{equation}\label{equ:Nwgd=0}
d\log|f_1|\wedge \cdots \wedge d\log |f_i| \wedge d\log f_{i+1} \wedge \cdots
\wedge d\log f_{2n} = 0.
\end{equation}
Denote by $x_j$ the left side of \eqref{equ:xj} multiplied by 
$\sqrt{-1}^{2n-2j-1}$. Taking the  imaginary part of 
\eqref{equ:Nwgd=0} and alternating $f_1,\cdots,f_{2n}$ we get
\begin{equation}\label{equ:sumxj=0}
\sum_{j}^{n-1}{2n-i \choose 2j+1-i} 
x_j=0\end{equation}
where the sum is over $j$ such that $2j+1 \geq i$.  
Indeed, denote the expression inside of ${\rm Alt}$ in \eqref{equ:xj} by 
$T_j$. An alternation of \eqref{equ:Nwgd=0} contributes to $T_j$ 
if and only if  $2j+1\ge i$, so that we can specify $2j+1-i$ terms from
$f_{i+1},\cdots,f_{2n}$ and make them contribute the $\log$ part of
$T_j$, which, together with $\log|f_1|, \cdots,\log|f_i|$, contribute
the log part in $T_j$. 

Let $s_i$ be the left hand side of \eqref{equ:sumxj=0} considered for
arbitrary $i$. Let us multiply it by ${2n\choose i} t^i$ and take a
sum over $0 \leq i \leq 2n-1$. Since $s_i =0$ for any
$0\le i\le n-1$ we have
$\sum_{i=0}^{2n-1} s_i{2n\choose i} t^i=t^n A(t)$ for some polynomial
$A(t)$ in $t$ whose coefficients are ${\Bbb Q}$-linear combinations of $x_i$'s. 
Using the identity
\begin{equation}\label{equ:comId}
{2n-i \choose p-i}{2n \choose i} = {2n \choose p}{p \choose i}.
\end{equation}
we have
$$\aligned
t^n A(t) =& \sum_{i=0}^{2n-1} \sum_{j}^{n-1}
{2n-i \choose 2j+1-i}{2n\choose i} x_jt^i  \\
=&\sum_{j=0}^{n-1}x_j{2n\choose 2j+1}\sum_{i=0}^{2j+1}{2j+1\choose i} t^i
=\sum_{j=0}^{n-1}x_j{2n\choose 2j+1}(1+t)^{2j+1}\text{\raisebox{-1ex}{.}}
\endaligned$$

Replacing $t$ by $t-1$ we get 
$\sum_{j=0}^{n-1}  x_j{2n\choose 2j+1} t^{2j+1} = (t-1)^n A(t-1)$. 
The left hand side is an odd polynomial of degree $2n-1$; it has 
a zero of order $n$ at $t=1$; therefore it must have a zero of order $n$
at $t= -1$, so its degree is at least $2n$. 
Thus it is a zero polynomial. 

(II) Denote the left side of \eqref{equ:yj} by $y_j$. Taking 
the real part of equation \eqref{equ:Nwgd=0} and alternating we get 
\begin{equation}\label{equ:sumyj=0}
\sum_{j=0}^n{2n-i\choose 2j-i}(-1)^{n-j}y_j = 0,
\qquad \qquad 0\le i\le n-1.
\end{equation}
By definition it is clear that $b_{n,n}= (2n)!$. 
Multiplying \eqref{equ:sumyj=0} by
${2n\choose i}t^i$ and taking sum over $0 \leq i \leq 2n$ we have
$$ \sum_{i=0}^{2n}\sum_{j\ge \lfloor i/2\rfloor}^n
(-1)^{n-j}y_j{2n-i\choose 2j-i}{2n\choose i} t^i
= t^nB(t)$$
for some polynomial $B(t)$. Using combinatorial identity \eqref{equ:comId}
we write it as 
$$
\sum_{j=0}^n (-1)^{n-j}y_j{2n\choose 2j}\sum_{i=0}^{2j}{2j\choose i} t^i 
= \sum_{j=0}^n (-1)^{n-j}y_j{2n\choose 2j}(t+1)^{2j}= t^nB(t).$$
Changing $t$ to $t-1$ and noticing that the left hand side 
is an even polynomial we get
$$\sum_{j=0}^n (-1)^{n-j}y_j{2n\choose 2j}t^{2j}= (t-1)^nB(t-1)
         = (t^2-1)^nC(t).$$
Therefore $C(t)=y_n$ is a constant. 
Thus we finally have
$y_j/y_n= b_{j,n}/(2n)!= {n\choose j}/{2n\choose 2j}$ ($0\le j\le n$).
The lemma is proved.
\end{proof}

Recall that the Grassmannian $n$-logarithm is defined by
$$
{\cal L}_n^{\rm G}(l_0,\dots,l_{2n-1})=(2\pi i)^{1-n}
\int_{{\Bbb CP}^{n-1}} r_{2n-1}(f_1,\dots, f_{2n-1})$$
where $f_i=l_i/l_0$.

\begin{prop} The Grassmannian $n$-logarithm 
${\cal L}_n^{\rm G}(l_0,\dots,l_{2n-1})$ can be expressed by
$$-\frac {(-2\pi i)^{1-n}}{(2n-1)!}
{\rm Alt}_{2n-1}\int_{{\Bbb CP}^{n-1}} {\rm Re}\left\{\log|f_1|
        \bigwedge_{j=2}^{n} d\log(f_j)\bigwedge_{j=n+1}^{2n-1}
        d\,\overline{\log(f_j)}\right\}$$
or
\begin{multline}\label{Grasslogn}
-\frac {(-4)^{n-1}((n-1)!)^2}{(2\pi i)^{n-1} (2n-2)!}
\int_{{\Bbb CP}^{n-1}}\log|f_1|\bigwedge_{j=2}^{2n-1} d\log|f_j| \\
=\frac {(-4)^{n-1}((n-1)!)^2}{(2\pi i)^{n-1} (2n-2)!(2n-1)!}{\rm Alt}_{2n}
\left\{\int_{{\Bbb CP}^{n-1}}\log|l_0|\bigwedge_{j=1}^{2n-2} d\log|l_j|\right\}.
\end{multline}
\end{prop}
\begin{proof} The first expression follows directly from 
Corollary \ref{cor:rn}. Now we prove the second. 
By definition \eqref{equ:rn}
\begin{multline}\label{rng}
r_{2n-1}(f_1,\cdots,f_{2n-1})= -\sum_{k=0}^{n-1} (-1)^{n-k-1} c_{k,2n-1} \\
{\rm Alt}_{2n-1}\Bigl\{\log|f_1|\, \bigwedge_{j=2}^{2k+1} d\log |f_j|
        \bigwedge_{j=2k+2}^{2n-1} d\arg f_j\Bigr\}
\end{multline}
where $c_{k,2n-1}={2n-1\choose 2k+1}/(2n-1)!$. Now let us look at the terms
in the expansion of \eqref{rng} which correspond to the term 
\begin{equation}\label{g1}
\log|f_1| d\log |f_2|\wedge\cdots \wedge d\log |f_{2n-1}|.
\end{equation}
By Lemma \ref{lem:xjyj}(II), each term inside the sum of \eqref{rng}
with $\log|f_1|$ as the first factor contributes to \eqref{g1} as many as 
$(-1)^{n-k}c_{k,2n-1}b_{k,n-1}$ times.
So the total contribution to  \eqref{g1} from \eqref{rng} is 
$$\aligned
d_n=&\sum_{k=0}^{n-1}(-1)^{n-k}c_{k,2n-1}b_{k,n-1}
=\sum_{k=0}^{n-1} \frac{(-1)^{n-k}}{2k+1}{n-1\choose k}\\
=&-\int_0^1 (t^2-1)^{n-1} \,dt
=(-1)^{n}\frac {\Gamma(n)\Gamma(\frac 12)}{2\cdot\Gamma(n+\frac 12)}
=(-1)^{n}\frac {2^{2n-2}((n-1)!)^2}{(2n-1)!}.
\endaligned$$
Therefore
\begin{multline}\label{middle} 
{\cal L}_n^{\rm G}(l_0,\cdots,l_{2n-1})=(2\pi i)^{1-n}
        \int_{{\Bbb CP}^{n-1}}r_{2n-1}(f_1,\cdots,f_{2n-1})\\
=-\frac {(-4)^{n-1}((n-1)!)^2}{(2\pi i)^{n-1}(2n-1)!}
        \sum_{i=1}^{2n-1}\sigma_{1i}\left(\int_{{\Bbb CP}^{n-1}}
          \log|f_1| \bigwedge_{j=2}^{2n-1} d\log|f_j| \right)
\end{multline}
where $\sigma_{11}=$id and for $i\ne 1$
$$\sigma_{1i}F(f_1,\dots,f_{2n-1})=
-F(f_i,f_2,\dots,f_{i-1},f_1,f_{i+1},\dots,f_{2n-2}).$$
Now we observe that for any $2\le i\le 2n-1$ we have
\begin{multline*}
\int_{{\Bbb CP}^{n-1}}\log|f_1| \bigwedge_{j=2}^{2n-1} d\log|f_j|-
\sigma_{1i}\left(\int_{{\Bbb CP}^{n-1}}\log|f_1| 
  \bigwedge_{j=2}^{2n-1} d\log|f_j|\right)\\
=(-1)^i \int_{{\Bbb CP}^{n-1}} d(\log |f_1|\log|f_i|)\wedge d\log|f_2| \wedge\cdots
\widehat{d\log|f_i|}\cdots \wedge d\log|f_{2n-2}|=0.
\end{multline*}
Therefore
$$\sum_{i=1}^{2n-1}\sigma_{1i} 
\left(\int_{{\Bbb CP}^{n-1}}\log|f_1|\bigwedge_{j=2}^{2n-1} d\log|f_j|\right)
=(2n-1)\int_{{\Bbb CP}^{n-1}}\log|f_1|\bigwedge_{j=2}^{2n-1} d\log|f_j|$$
which together with equation \eqref{middle} 
yields the second equality. To prove the last equality in our
proposition it suffices to observe that
$$\aligned
(2n-1)!\int_{{\Bbb CP}^{n-1}} \log|f_1|\bigwedge_{j=2}^{2n-1}d\log|f_j|
=&{\rm Alt}_{2n-1}\int_{{\Bbb CP}^{n-1}} 
        \log|f_1|\bigwedge_{j=2}^{2n-1} d\log|f_j|\\
=&-{\rm Alt}_{2n}  \int_{{\Bbb CP}^{n-1}} 
        \log|l_0|\bigwedge_{j=1}^{2n-2} d\log|l_j|.
\endaligned$$
\end{proof}

\begin{rem}
This result improves Proposition 3.2 of \cite{Gchow}.
\end{rem}

\end{document}